\documentclass{article}

\usepackage{amsmath,amssymb,amsthm}
\usepackage{mathrsfs}
\usepackage{amscd}
\usepackage{enumerate}

\theoremstyle{plain}
\newtheorem{theorem}{Theorem}[section]
\theoremstyle{remark}

\theoremstyle{plain}
\newtheorem{corollary}[theorem]{Corollary}

\numberwithin{equation}{section}

\def\R{{\mathbb R}}

\newcommand{\uh}{\underline{h}}
\newcommand{\oh}{\overline{h}}

\usepackage{color}

\definecolor{gr}{rgb}   {0.,   0.8,   0. } 
\definecolor{bl}{rgb}   {0.,   0.5,   1. } 
\definecolor{mg}{rgb}   {0.7,  0.,    0.7}

\begin{document}

\title{Change of angles in tent spaces}

\date{January 14, 2011}

\author{Pascal Auscher}

\maketitle

\begin{abstract}

We prove sharp bounds for the equivalence of norms in tent spaces with respect to changes of angles. Some applications are given.

\end{abstract}

\bigskip

Let $B(x,t)$ denote an open ball centered at $x\in \R^n$ with radius $t>0$.
Define for a locally square integrable function $g(t,y)$, $(t,y)\in \R^{n+1}_{+}$,  for  $\alpha>0$ and $x\in \R^n$,
\begin{equation}
\label{eq:a}
A^{(\alpha)}g(x) := \bigg(  
\int_{\R^{n+1}_{+}} \frac{1_{B(x,\alpha t)}(y)}{t^{n}}\, \bigl|g(t,y)\bigr|^{2} 
\,\frac{dy dt}{t}\bigg)^{\frac{1}{2}},
\end{equation}
 and for $0<p<\infty$, say that $g\in T^{p,2}_\alpha$ if 
$$\|g\|_{T^{p,2}_\alpha}:= \|A^{(\alpha)}g\|_{p}<\infty.$$
This space was introduced in \cite{cms}. As sets the spaces for a given $p$ are the same and  the norms (or quasi-norms if $p<1$) are equivalent: whenever  $\alpha,\beta>0$, one has
 \begin{equation}
\label{eq:equiv}
\|A^{(\alpha)}g\|_{p} \sim  \|A^{(\beta)}g\|_{p}\ .
\end{equation}

For $p=\infty$, the limiting space is defined with a Carleson measure condition. Let  $T_{\alpha}B$ be the tent with aperture $\alpha$ above the open ball $B=B(x,r)$, \textit{i.e.,\, }the set of $(t,y)$ such that 
 $0<t< r/\alpha$ and $y\in B(x, r-\alpha t)$.  
 We define $\|g\|_{T^{\infty,2}_\alpha}$ as the infimum of  $C\ge 0$ such that for all ball $B$,
 \begin{equation}
\label{eq:b}
\int_{T_{\alpha}B} |g(t,y)|^{2}\, \frac{dy dt}{t} \le C^2\, \frac{|B|}{\alpha^n}\ .
\end{equation}
 Again, the spaces $T^{\infty,2}_{\alpha}$ are the same, the norms are equivalent and the isometry property holds. 
 
 To expain the choice of the normalisation in (\ref{eq:b}), we remark that for $p=\infty$ included,   
  $g(t,y)\mapsto h(t,y):=\alpha^{n/2} g(t/\alpha,y)$ is an isometry between $T^{p,2}_{1}$ and $T^{p,2}_{\alpha}$ equipped with their respective (quasi-)norms. 
 
 It is shown in \cite{cms} that the spaces $T^{p,2}_{1}$, $0<p\le\infty$, interpolate by the complex method and the real method. The same results hold for $T^{p,2}_\alpha$ for fixed $\alpha$, with constants (\textit{i.e.,\,}the constants in the equivalence of norms between $T^{p,2}_\alpha$ and the interpolated space to which it is equal) independent of $\alpha$ by using the isometry property. 
 
 Motivated by an intensive usage of tent spaces in the  development of new Hardy spaces associated to operators with Gaffney-Davies estimates first made in \cite{amr}, \cite{HM}, and also by the study of maximal regularity on tent spaces towards applications for parabolic PDE's (\cite{AMP10}, and some more work in progress), it became interesting to know the sharp dependence of the bounds in (\ref{eq:equiv}) with respect to $\alpha,\beta$. The $L^2$ bound is immediate by Fubini's theorem:  $\|g\|_{T^{2,2}_\alpha}=(\alpha/\beta)^{n/2}\|g\|_{T^{2,2}_{\beta}}$. For $p\ne 2$, 
 the argument in \cite{cms}, originating from  \cite{fs},  does not give optimal dependence. 
 The inequality 
 $$  
    \| g\|_{T_{\alpha}^{p,2}}\leq C(1+\log\alpha)\alpha^{n/\tau}
   \| g\|_{T^{p,2}_{1}}\ ,
$$
for $1<p<\infty$ and $\alpha \geq 1$, where $\tau = \min(p,2)$ and $C$ depends only on $n$ and $p$,  is proved in \cite{hnp} as a special case of a Banach space valued result, and, moreover, the polynomial growth $\alpha^{n/\tau}$ is shown to be optimal for such an inequality to hold.   
The restriction $p>1$ occurs  in this argument because the UMD property is required  and a maximal inequality is used. Note that even the $L^2$ bounds is not immediate in a Banach (non Hilbert) space valued context. In discussion with T. Hyt\"onen, we convinced ourselves that the logarithmic factor is not produced by this argument in the scalar case when $p\ge 2$.  Still, this argument is quite involved and elimination of the logarithm in the $p<2$ situation was unclear.

Here we give the sharp lower and upper bounds for (\ref{eq:equiv}) in the scalar case by a very simple argument.  Define $\uh(p,\alpha)=\min\big( \alpha^{-n/2}, \alpha^{-n/p}\big)$, $\oh(p,\alpha)=\max\big( \alpha^{-n/2}, \alpha^{-n/p}\big)$. Note that $\uh(p,\alpha)= \alpha^{-n/p}$ if $(\alpha-1)(p-2)\ge 0$ and $\uh(p,\alpha)= \alpha^{-n/2}$ if $(\alpha-1)(p-2)\le 0$, and inversely for $\oh(p,\alpha)$.

\begin{theorem}\label{thm:anglebetter}
Let $0<p\le\infty$ and $\alpha,\beta>0$. There exist  constants $C,C'>0$ depending on $n,p$ only,  such that for any locally square integrable function $g$, 
$$
C\uh(p,\alpha/\beta)   \| g\|_{T^{p,2}_{\alpha}}\leq \| g\|_{T_{\beta}^{p,2}}\leq C'\oh(p,\alpha/\beta)   \| g\|_{T^{p,2}_{\alpha}}\ .
$$
Moreover, the dependence in $\alpha/\beta$ is best possible in the sense that this growth is attained. \end{theorem}

In particular,  for $\alpha>1$, $\tau=\min(2,p)$ and $\sigma=\max(2,p)$, one has \begin{equation}
\label{eq:1}
\| g\|_{T_{\alpha}^{p,2}}\leq C\alpha^{n/\tau}
   \| g\|_{T^{p,2}_{1}}\ ,
\end{equation}
 \begin{equation}
\label{eq:2}
 \| g\|_{T_{1}^{p,2}}\leq C\alpha^{-n/\sigma}
   \| g\|_{T^{p,2}_{\alpha}}\  . 
\end{equation} The second one improves  the obvious bound $\| g\|_{T_{1}^{p,2}}\leq 
   \| g\|_{T^{p,2}_{\alpha}}$. By symmetry using the relation $\uh(p,\alpha)^{-1}= \oh(p, \alpha^{-1})$ and scale invariance, all cases reduce to (\ref{eq:1}) and (\ref{eq:2}) with $\alpha>1$.  

\begin{corollary}\label{thm:VS}
Let $0<p<\infty$. There is  a constant $0<C<\infty$ depending on $n,p$ only such that for any locally square integrable $g$, if $Vg(x)=\left(\int_{0}^\infty |g(t,x)|^2\, \frac{dt}{t}\right)^{\frac{1}{2}},$ %
$$
\|Vg\|_{p}\leq C   \| g\|_{T^{p,2}_{1}}, \quad \mathrm{if} \ p\le 2, \qquad 
 \|Vg\|_{p} \geq C 
   \| g\|_{T^{p,2}_{1}},  \quad \mathrm{if} \ p\ge 2.
$$
\end{corollary}

The corollary was proved by a different method in \cite{AHM10} when $p<2$ and the $p>2$ case dates back to \cite{St}.  The opposite inequalities are false if $p\ne 2$. Starting from Theorem \ref{thm:anglebetter}, the proof is a mere application of Lebesgue differentiation theorem when $\alpha\to 0$ for $\alpha^{-n/2}A^{(\alpha)}g(x)$ converges to $Vg(x)$  assuming $g$ smooth with compact support. This assumption is easily removed.

\begin{corollary} If $0<p<\infty$ and  $\lambda>\max(2/p,1)$, for any locally square integrable $g$, 
$$
\bigg\| \biggl( 
\int_{\R^{n+1}_{+}} \left(\frac{t}{|x-y|+t}\right)^{n\lambda}\, \bigl|g(t,y)\bigr|^{2} 
\,\frac{dy dt}{t^{n+1}}\biggr)^{\frac{1}{2}}\bigg\|_{p} \le C(n,p, \lambda) \|A^{(1)}g\|_{p}\ .
$$

\end{corollary}

The left hand side equals the grand square function of Stein when    $g(t,y)=t\nabla u(t,y)$, $u$ being  the harmonic extension of a suitable  distribution $u_{0}$ on $\R^n$.  Hence,  for all $0<p<\infty$ and  $\lambda>\max(2/p,1)$,  it is dominated in $L^p$ by  $ \|A^{(1)}g\|_{p}$ which is the $L^p$ norm of  the area functional of Lusin defined from $u_{0}$. However, it is known from Stein-Weiss' theory that   $A^{(1)}g \in L^p(\R^n)$ if and only if  $\frac{n-1}n<p<\infty$ and $u_{0}$ belongs to the Hardy space $H^p(\R^n)$  (See \cite{St}).     
This gives a simple proof of Theorem 2, Chap. IV in \cite{St}. The lower exponent $\frac{n-1}{n}$ is only due to the choice of the extension. Using an extension by convolution with  $t^{-n}\varphi(x/t)$ with $\varphi\in C^\infty_{0}(\R^n)$ having all vanishing moments but the one of order 0, $\frac{n-1}{n}$ becomes 0 by the results in \cite{fs}. At $\lambda=2/p$ and $p<2$,  a weak type inequality  is plausible,   the Lorentz norm $L^{p,\infty}$ replacing the Lebesgue norm $L^p$ in  the left hand side. It would give a simple proof of  the weak type $(p,p)$ result of Fefferman \cite{Fef} for Stein's grand square function.  We leave this open. 

The proof of the corollary is easy by splitting the upper half space according to  $|x-y|/t$ compared to powers  $2^k$, $k\ge 0$, and one obtains
$$
\biggl(  
\int_{\R^{n+1}_{+}} \left(\frac{t}{|x-y|+t}\right)^{n\lambda}\, \bigl|g(t,y)\bigr|^{2} 
\,\frac{dy dt}{t^{n+1}}\biggr)^{\frac{1}{2}} \le C(n,\lambda)\sum_{k\ge 0} 2^{-kn\lambda/2} A^{(2^{k+1})}g(x)\ .$$
It remains to use $    \| A^{(\alpha)}g\|_{p}\leq C\alpha^{n/\tau}   \| A^{(1)}g\|_{p}
$ for $\alpha=2^{k+1}$ in appropriate arguments for $p\ge 1$ or $p\le 1$.

The proof of Theorem \ref{thm:anglebetter} is an easy matter  using in part atomic theory for tent spaces, again proved in \cite{cms}. Recall that for $0<p\le 1$, the tent space $T^{p,2}_{1}$ has an atomic decomposition: A $T^{p,2}_{1}$ atom is a function $a(t,x)$ supported in a tent $T_{1}B$   with the estimate
$$
 \int_{T_{1}B} |a(t,x)|^{2}\,  \frac{dx dt}{t} \le |B|^{-(\frac{2}{p}-1)}\ .
 $$
 There is a constant $C=C(n,p)>0$ such that $\|a\|_{T^{p,2}_{1}}\le C$.
 Any $T^{p,2}_{1}$ function $g$ can be represented as a  series $g=\sum \lambda_{j }a_{j}$ where $a_{j}$ is a $T^{p,2}_{1}$ atom and $\sum |\lambda_{j}|^p \sim \|g\|_{T^{p,2}_{1}}^p$. By the isometry property, a $T^{p,2}_{\alpha}$ atom is a function   $A(u,x)$  supported  a tent $T_{\alpha}B$   with the estimate
$$
 \int_{T_{\alpha}B}  |A(u,x)|^{2} \, \frac{ dx du}{u} \le \alpha^{-n}|B|^{-(\frac{2}{p}-1)}
 $$
 and the decomposition theorem holds in $T^{p,2}_{\alpha}$.   
%
%

{\bf Proof of Theorem \ref{thm:anglebetter}.} Recall that we may assume $\alpha>\beta=1$ and it is enough to  prove (\ref{eq:1}) and (\ref{eq:2}).     
Fix $p=\infty$ first.  Let $\|g\|_{T^{\infty,2}_{1}}=1$ and $B$ be a ball. As $T_{\alpha}B \subset T_{1}B$\ , 
$$
 \int_{T_{\alpha}B} |g(t,x)|^{2}\, \frac{dx dt}{t}  \le  \int_{T_{1}B} |g(t,x)|^{2} \, \frac{dx dt}{t} \le |B| =  \alpha^n \frac{|B|}{\alpha^n}\ .
 $$
Hence, $\|g\|_{T^{\infty,2}_\alpha} \le \alpha^{n/2}$. This shows $\|g\|_{T^{\infty,2}_\alpha} \le \alpha^{n/2} \|g\|_{T^{\infty,2}_1} $ for all $g$.

Let $\|g\|_{T^{\infty,2}_{\alpha}}=1$ and $B$ be a ball. As $T_{1}B \subset T_{\alpha}(\alpha B) $ where $\alpha B$ is the ball concentric with $B$ dilated by $\alpha$, 
$$
 \int_{T_{1}B} |g(t,x)|^{2}\,  \frac{dx dt}{t}  \le  \int_{T_{\alpha}(\alpha B)} |g(t,x)|^{2}\, \frac{dx dt}{t} \le   \frac{|\alpha B|}{\alpha^n}= |B|\ .
 $$
 Hence, $\|g\|_{T^{\infty,2}_1} \le1$ and 
$\|g\|_{T^{\infty,2}_1} \le  \|g\|_{T^{\infty,2}_\alpha} $ for all $g$.

Fix now $p\le 1$. Let  $B$ be a	 ball and $a$ be a $T^{p,2}_{1}$ atom supported in a tent $T_{1}B$.  As $T_{1}B \subset T_{\alpha}(\alpha B)$, we have $a$ is supported in $T_{\alpha}(\alpha B)$ and 
$$
 \int_{T_{\alpha}(\alpha B) }  |a(t,x)|^{2}\,  \frac{ dx dt}{t} = \int_{T_{1}B} |a(t,x)|^{2}\,  \frac{dx dt}{t} \le |B|^{-(\frac{2}{p}-1)}= \alpha^{2n/p} \big(\alpha^{-n}|\alpha B|^{-(\frac{2}{p}-1)}\big)\ .
 $$
 Thus $\alpha^{-n/p} a$ is a $T^{p,2}_{\alpha}$ atom. An atomic decomposition of any element of $T^{p,2}_{1}$ is up to multiplication by $\alpha^{-n/p}$  an atomic decomposition in $T^{p,2}_{\alpha}$,  proving $\|g\|_{T^{p,2}_\alpha} \le C(n,p)\alpha^{n/p} \|g\|_{T^{p,2}_1} $.
 
Next, let  $B$ be a	 ball and $a$ be a $T^{p,2}_{\alpha}$ atom supported in a tent $T_{\alpha}B$.  As $T_{\alpha}B \subset T_{1}B$, we have $a$ is supported in $T_{1}B$ and 
$$
 \int_{T_1B }  |a(t,x)|^{2}\,  \frac{ dx dt}{t} = \int_{T_{\alpha} B} |a(t,x)|^{2}\,  \frac{dx dt}{t} \le \alpha^{-n}  |B|^{-(\frac{2}{p}-1)}\ .
 $$
 Thus $\alpha^{n/2} a$ is a $T^{p,2}_{1}$ atom. As above, we conclude that $\|g\|_{T^{p,2}_1} \le C(n,p)\alpha^{-n/2} \|g\|_{T^{p,2}_\alpha} $.

 For $1< p< 2$ and $2<p<\infty$, we conclude by interpolation with the $p=2$ equality   $\|g\|_{T^{2,2}_\alpha}=\alpha^{n/2}\|g\|_{T^{2,2}_{1}}$. We have shown (\ref{eq:1}) and (\ref{eq:2}).

 The sharpness of the bounds is seen by saturating these inequalities. Fix $\alpha>1$ large. Let $B$ be the unit ball. Set $a_{1}(t,y)=1_{T_{1}B}(t,y) 1_{[1/2,1]}(t).$ It is easy to see that $\|a_{1}\|_{T^{p,2}_{1}}\sim 1$. Now, we have that $A^{(\alpha)}a_{1}$ has support equal to $\overline{B}(0,\alpha)$, is bounded by a constant $c(n)>0$ and equal  to that  constant  on the ball $B(0, \frac{\alpha+1}{2})$. Thus $\|a_{1}\|_{T^{p,2}_{\alpha}}\sim \alpha^{n/p}$. This proves that (\ref{eq:1})  is optimal when $p\le 2$ and (\ref{eq:2}) is optimal when $p\ge 2 $. 
 Next, let $a_{2}(t,y)= a_{1}(\alpha t, y)$.  By scaling $\|a_{2}\|_{T^{p,2}_{\alpha}}=\alpha^{n/2}\|a_{1}\|_{T^{p,2}_{1}} \sim \alpha^{n/2}$. Now, $A^{(1)}a_{2}$ is supported in $\overline B$, bounded by a constant $c(n)>0$ and equal to that constant on 
 $B(0, \frac{\alpha-1}{2\alpha})$. Thus $\|a_{2}\|_{T^{p,2}_{1}}  \sim 1$. This proves that (\ref{eq:1})  is optimal when $p\ge 2$ and (\ref{eq:2}) is optimal when $p\le 2 $. 
\qed




\end{document}